\documentclass[11pt]{article}
\usepackage{amssymb,amsthm,amsmath}
\usepackage[dvipdfm,colorlinks,breaklinks,
pdfstartview={FitH -32768},
bookmarks,bookmarksnumbered,bookmarksopen,
pdftitle={Classification of outer actions of Z^N on O_2},
pdfkeywords={2000 Mathematics Subject Classification: 46L55.},
pdfauthor={Hiroki Matui}]{hyperref}

\newtheorem{thm}{Theorem}[section]
\newtheorem{prop}[thm]{Proposition}
\newtheorem{lem}[thm]{Lemma}
\newtheorem{cor}[thm]{Corollary}

\theoremstyle{definition}
\newtheorem{df}[thm]{Definition}
\newtheorem{rem}[thm]{Remark}

\numberwithin{equation}{section}

\newcommand{\N}{\mathbb{N}}
\newcommand{\Z}{\mathbb{Z}}
\newcommand{\R}{\mathbb{R}}
\newcommand{\C}{\mathbb{C}}

\newcommand{\Aut}{\operatorname{Aut}}
\newcommand{\Ad}{\operatorname{Ad}}
\newcommand{\Sp}{\operatorname{Sp}}

\newcommand{\ep}{\varepsilon}

\title{Classification of outer actions of $\Z^N$ on $\mathcal{O}_2$}
\author{Hiroki Matui 
\thanks{Supported in part by a grant 
from the Japan Society for the Promotion of Science} \\
Graduate School of Science \\
Chiba University \\
1-33 Yayoi-cho, Inage-ku, Chiba 263-8522, Japan}
\date{}

\begin{document}
\maketitle

\begin{abstract}
We will show that 
any two outer actions of $\Z^N$ on $\mathcal{O}_2$ 
are cocycle conjugate. 
\end{abstract}

\section{Introduction}

Group actions on $C^*$-algebras and von Neumann algebras are 
one of the most fundamental subjects in the theory of operator algebras. 
A. Connes introduced a non-commutative Rohlin property 
and classified single automorphisms of von Neumann algebras 
(\cite{C1},\cite{C2}), 
and A. Ocneanu generalized it 
to actions of discrete amenable groups (\cite{O}). 
In the setting of $C^*$-algebras, 
A. Kishimoto established a non-commutative Rohlin type theorem 
for single automorphisms on UHF algebras 
and classified them up to outer conjugacy (\cite{K1},\cite{K2}). 
H. Nakamura used the same technique 
for automorphisms on purely infinite simple nuclear $C^*$-algebras 
and classified them up to outer conjugacy (\cite{N2}). 
Furthermore, he proved a Rohlin type theorem 
for $\Z^2$-actions on UHF algebras 
and gave a classification result for product type actions 
up to cocycle conjugacy (\cite{N1}). 
Recently, T. Katsura and the author generalized this result 
and gave a complete classification of uniformly outer actions 
of $\Z^2$ on UHF algebras (\cite{KM}). 
In the case of finite group actions, 
M. Izumi introduced a notion of the Rohlin property 
and classified a large class of actions (\cite{I2},\cite{I3}). 
The reader may consult the survey paper \cite{I1} 
for the Rohlin property of automorphisms on $C^*$-algebras. 

The aim of this paper is 
to extend these results to $\Z^N$-actions 
on the Cuntz algebra $\mathcal{O}_2$. 
More precisely, we will show that 
any outer actions of $\Z^N$ on $\mathcal{O}_2$ have the Rohlin property 
and that they are cocycle conjugate to each other. 
The content of this paper is as follows. 
In Section 2, 
we collect notations and basic facts relevant to this paper. 
Notions of the ultraproduct algebra $A^\omega$ and 
the central sequence algebra $A_\omega$ will help our analysis. 
When $A$ is isomorphic to $\mathcal{O}_2$, 
it is known that $A_\omega$ contains a unital copy of $\mathcal{O}_2$. 
This fact implies 
strong triviality of the $K$-theory of $\mathcal{O}_2$. 
For $\Z^N$-actions on unital $C^*$-algebras, 
we recall the definition of the Rohlin property from \cite{N1} 
and give a couple of remarks. 

In Section 3, we establish the cohomology vanishing theorem 
for $\Z^N$-actions on the Cuntz algebra $\mathcal{O}_2$ 
with the Rohlin property. 
One of the difficulties 
in the study of $\Z^N$-actions on $C^*$-algebras is that 
one has to deal with homotopies of unitaries 
in order to obtain the so-called cohomology vanishing theorem. 
In our situation, however, 
we do not need 
to take care of $K_*$-classes of (continuous families of) unitaries, 
because of the triviality of $K_*(\mathcal{O}_2)$. 
Indeed, any continuous family of unitaries in $\mathcal{O}_2$ 
is homotopic to the identity by a smooth path of finite length. 
This enables us to avoid $K$-theoretical arguments. 

Section 4 is devoted to the Rohlin type theorem 
for $\Z^N$-actions on $\mathcal{O}_2$. 
The main idea of the proof is similar to that in \cite{N1}. 
We also make use of several techniques developed in \cite{N2}. 
The proof is by induction on $N$, 
because we need the cohomology vanishing theorem for $\Z^{N-1}$-actions 
in order to prove the Rohlin type theorem for $\Z^N$-actions. 

In Section 5, the main theorem is shown. 
D. E. Evans and A. Kishimoto introduced in \cite{EK} 
an intertwining argument for automorphisms, 
which is an equivariant version of Elliott's intertwining argument 
for classification of $C^*$-algebras. 
By using the Evans-Kishimoto intertwining argument, 
we show that 
any two outer actions of $\Z^N$ on $\mathcal{O}_2$ are cocycle conjugate. 
The cohomology vanishing theorem is necessary 
in each step of the intertwining argument. 

\bigskip

\textbf{Acknowledgement.} 
The author is grateful to Toshihiko Masuda for many helpful discussions.

\section{Preliminaries}

Let $A$ be a unital $C^*$-algebra. 
We denote by $\mathcal{U}(A)$ the group of unitaries in $A$. 
For $u\in\mathcal{U}(A)$, we let $\Ad u(a)=uau^*$ for $a\in A$ and 
call it an inner automorphism on $A$. 
When an automorphism $\alpha\in\Aut(A)$ on $A$ is not inner, 
it is said to be outer. 

For any $a,b\in A$, we write $[a,b]=ab-ba$ 
and call it the commutator of $a$ and $b$. 

\bigskip

In this paper, we deal with central sequence algebras, 
which simplify the arguments a little. 
Let $A$ be a separable $C^*$-algebra and 
let $\omega\in\beta\N\setminus\N$ be a free ultrafilter. 
We set 
\[ c^\omega(A)=\{(a_n)\in\ell^\infty(\N,A)\mid
\lim_{n\to\omega}\lVert a_n\rVert=0\}, \]
\[ A^\omega=\ell^\infty(\N,A)/c^\omega(A). \]
We identify $A$ with the $C^*$-subalgebra of $A^\omega$ 
consisting of equivalence classes of constant sequences. 
We let 
\[ A_\omega=A^\omega\cap A'. \]
When $\alpha$ is an automorphism on $A$ or 
an action of a discrete group on $A$, 
we can consider its natural extension on $A^\omega$ and $A_\omega$. 
We denote it by the same symbol $\alpha$. 

If $A$ is a unital separable purely infinite simple $C^*$-algebra, 
then $A^\omega$ is purely infinite simple. 
When $A$ is 
a unital separable purely infinite simple nuclear $C^*$-algebra, 
it is known that $A_\omega$ is also purely infinite simple 
(\cite[Proposition 3.4]{KP}). 

\bigskip

We would like to collect several facts 
about the Cuntz algebra $\mathcal{O}_2$. 
The Cuntz algebra $\mathcal{O}_2$ is 
the universal $C^*$-algebra generated by two isometries $s_1$ and $s_2$ 
satisfying $s_1s_1^*+s_2s_2^*=1$. 
It is a unital separable purely infinite simple nuclear $C^*$-algebra 
with trivial $K$-groups, i.e. $K_0(\mathcal{O}_2)=K_1(\mathcal{O}_2)=0$. 
By \cite[Theorem 3.6]{R1}, 
any automorphisms $\alpha$ and $\beta$ on $\mathcal{O}_2$ 
are approximately unitarily equivalent. 
Thus, 
there exists a sequence of unitaries $\{u_n\}_n$ in $\mathcal{O}_2$ 
such that $\beta(a)=\lim_{n\to\infty}u_n\alpha(a)u_n^*$ 
for all $a\in A$. 
It is also known that 
$\mathcal{O}_2$ is isomorphic to the infinite tensor product 
$\bigotimes_{n=1}^\infty\mathcal{O}_2$ 
(\cite{R2} or \cite[Theorem 3.8]{KP}). 
Consequently, 
$(\mathcal{O}_2)_\omega=(\mathcal{O}_2)^\omega\cap\mathcal{O}_2'$ 
contains a unital copy of $\mathcal{O}_2$. 
(Actually, its converse also holds: 
if $A$ is a unital simple separable nuclear $C^*$-algebra 
and $A_\omega$ contains a unital copy of $\mathcal{O}_2$, 
then $A$ is isomorphic to $\mathcal{O}_2$ (\cite[Lemma 3.7]{KP}). 
But, in this paper, we do not need this fact.) 
Let $e$ be a projection of $(\mathcal{O}_2)_\omega$. 
By the usual argument taking subsequences, 
we can find a unital copy of $\mathcal{O}_2$ 
in the relative commutant of $e$ in $(\mathcal{O}_2)_\omega$. 
Therefore, $[e]=[e]+[e]$ in $K_0((\mathcal{O}_2)_\omega)$, 
and so $[e]=0$. 
Since $e$ is arbitrary and 
$(\mathcal{O}_2)_\omega$ is purely infinite simple, 
we have $K_0((\mathcal{O}_2)_\omega)=0$. 
In a similar fashion, we also have $K_1((\mathcal{O}_2)_\omega)=0$. 
Indeed, for any unitary $u\in(\mathcal{O}_2)_\omega$, 
we can find a unital copy of $\mathcal{O}_2$ 
in the relative commutant of $u$ in $(\mathcal{O}_2)_\omega$. 
It follows from \cite[Lemma 5.1]{HR} that 
there exists a smooth path $u(t)$ of unitaries 
in $(\mathcal{O}_2)_\omega$ such that $u(0)=u$ and $u(1)=1$. 
Moreover, the length of the path is not greater than $8\pi/3$. 
We will use this argument in Lemma \ref{HaagerupRordam} again. 

\bigskip

Let $N$ be a natural number and 
let $\xi_1,\xi_2,\dots,\xi_N$ be the canonical basis of $\Z^N$, 
that is, 
\[ \xi_i=(0,0,\dots,1,\dots,0,0), \]
where $1$ is in the $i$-th component. 
Let $\alpha$ be an action of $\Z^N$ on a unital $C^*$-algebra $A$. 
We say that $\alpha$ is an outer action on $A$, 
if $\alpha_n$ is not inner on $A$ for any $n\in\Z^N\setminus\{0\}$. 

A family of unitaries $\{u_n\}_{n\in\Z^N}$ in $A$ is called 
an $\alpha$-cocycle, 
if 
\[ u_n\alpha_n(u_m)=u_{n+m} \]
for all $n,m\in\Z^N$. 
If a family of unitaries $u_1,u_2,\dots,u_N\in\mathcal{U}(A)$ satisfies 
\[ u_i\alpha_{\xi_i}(u_j)=u_j\alpha_{\xi_j}(u_i) \]
for all $i,j=1,2,\dots,N$, 
then it determines uniquely an $\alpha$-cocycle $\{u_n\}_{n\in\Z^N}$ 
such that $u_{\xi_i}=u_i$. 
We may also call the family $\{u_1,u_2,\dots,u_N\}$ an $\alpha$-cocycle. 
An $\alpha$-cocycle $\{u_n\}_n$ in $A$ is called a coboundary, 
if there exists $v\in\mathcal{U}(A)$ such that 
\[ u_n=v\alpha_n(v^*) \]
for all $n\in\Z^N$, or equivalently, 
if 
\[ u_{\xi_i}=v\alpha_{\xi_i}(v^*) \]
for all $i=1,2,\dots,N$. 

When $\{u_n\}_{n\in\Z^N}$ is an $\alpha$-cocycle, 
it turns out that 
a new action $\tilde{\alpha}$ of $\Z^N$ on $A$ can be defined by 
\[ \tilde{\alpha}_n(x)=\Ad u_n\circ\alpha_n(x)=u_n\alpha_n(x)u_n^* \]
for each $x\in A$. 
We call $\tilde{\alpha}$ 
the perturbed action of $\alpha$ by $\{u_n\}_n$. 

Two actions $\alpha$ and $\beta$ of $\Z^N$ on $A$ 
are said to be cocycle conjugate, 
if there exists an $\alpha$-cocycle $\{u_n\}_n$ in $A$ such that 
the perturbed action of $\alpha$ by $\{u_n\}_n$ is conjugate to $\beta$. 
The main theorem of this paper states that 
any two outer actions of $\Z^N$ on $\mathcal{O}_2$ 
are cocycle conjugate to each other. 

\bigskip

Let $N$ be a natural number. 
We would like to recall 
the definition of the Rohlin property for $\Z^N$-actions 
on unital $C^*$-algebras (see \cite[Section 2]{N1}). 
Let $\xi_1,\xi_2,\dots,\xi_N$ be the canonical basis of $\Z^N$ as above. 
For $m=(m_1,m_2,\dots,m_N)$ and $n=(n_1,n_2,\dots,n_N)$ in $\Z^N$, 
$m\leq n$ means $m_i\leq n_i$ for all $i=1,2,\dots,N$. 
For $m=(m_1,m_2,\dots,m_N)\in\N^N$, we let 
\[ m\Z^N=\{(m_1n_1,m_2n_2,\dots,m_Nn_N)\in\Z^N
\mid (n_1,n_2,\dots,n_N)\in\Z^N\}. \]
For simplicity, we denote $\Z^N/m\Z^N$ by $\Z_m$. 
Moreover, we may identify $\Z_m=\Z^N/m\Z^N$ with 
\[ \{(n_1,n_2,\dots,n_N)\in\Z^N
\mid 0\leq n_i\leq m_i-1\text{ for all }i=1,2,\dots,N\}. \]

\begin{df}\label{Rohlin}
Let $\alpha$ be an action of $\Z^N$ on a unital $C^*$-algebra $A$. 
Then $\alpha$ is said to have the Rohlin property, 
if for any $m\in\N$ there exist $R\in\N$ and 
$m^{(1)},m^{(2)},\dots,m^{(R)}\in\N^N$ 
with $m^{(1)},\dots,m^{(R)}\geq(m,m,\dots,m)$ 
satisfying the following: 
For any finite subset $\mathcal{F}$ of $A$ and $\ep>0$, 
there exists a family of projections 
\[ e^{(r)}_g \qquad (r=1,2,\dots,R, \ g\in\Z_{m^{(r)}}) \]
in $A$ such that 
\[ \sum_{r=1}^R\sum_{g\in\Z_{m^{(r)}}}e^{(r)}_g=1, 
\quad \lVert[a,e^{(r)}_g]\rVert<\ep, 
\quad \lVert\alpha_{\xi_i}(e^{(r)}_g)-e^{(r)}_{g+\xi_i}\rVert<\ep \]
for any $a\in\mathcal{F}$, $r=1,2,\dots,R$, $i=1,2,\dots,N$ 
and $g\in\Z_{m^{(r)}}$, 
where $g+\xi_i$ is understood modulo $m^{(r)}\Z^N$. 
\end{df}

\begin{rem}\label{reindex}
Clearly, we can restate the definition of the Rohlin property 
as follows. 
For any $m\in\N$ there exist $R\in\N$, 
$m^{(1)},m^{(2)},\dots,m^{(R)}\in\N^N$ 
with $m^{(1)},\dots,m^{(R)}\geq(m,m,\dots,m)$ and 
a family of projections 
\[ e^{(r)}_g \qquad (r=1,2,\dots,R, \ g\in\Z_{m^{(r)}}) \]
in $A_\omega=A^\omega\cap A'$ such that 
\[ \sum_{r=1}^R\sum_{g\in\Z_{m^{(r)}}}e^{(r)}_g=1, 
\quad \alpha_{\xi_i}(e^{(r)}_g)=e^{(r)}_{g+\xi_i} \]
for any $r=1,2,\dots,R$, $i=1,2,\dots,N$ and $g\in\Z_{m^{(r)}}$, 
where $g+\xi_i$ is understood modulo $m^{(r)}\Z^N$. 

Furthermore, by the reindexation trick, 
for a given separable subset $S$ of $A^\omega$, 
we can make the projections $e^{(r)}_g$ commute 
with all elements in $S$. 
We refer the reader to \cite[Lemma 5.3]{O} for details. 

In particular, the same conclusion also follows 
for perturbed actions on $A^\omega$. 
Let $\alpha$ be an action of $\Z^N$ on $A$ with the Rohlin property 
and let $\{u_n\}_n\subset\mathcal{U}(A^\omega)$ be 
an $\alpha$-cocycle in $A^\omega$. 
We can consider 
the perturbed action $\tilde{\alpha}$ of $\Z^N$ on $A^\omega$. 
Then, for a given separable subset $S$ of $A^\omega$, 
we can choose the projections $e^{(r)}_g$ in $A_\omega$ 
so that they commute with all elements in $S\cup\{u_n\}_n$. 
It follows that we have 
\[ \tilde{\alpha}_{\xi_i}(e^{(r)}_g)
=u_{\xi_i}\alpha_{\xi_i}(e^{(r)}_g)u_{\xi_i}^*
=u_{\xi_i}e^{(r)}_{g+\xi_i}u_{\xi_i}^*
=e^{(r)}_{g+\xi_i} \]
for any $r=1,2,\dots,R$, $i=1,2,\dots,N$ and $g\in\Z_{m^{(r)}}$, 
where $g+\xi_i$ is understood modulo $m^{(r)}\Z^N$. 
\end{rem}

\begin{rem}\label{Rem2ofN1}
We can also restate the Rohlin property as follows (\cite[Remark 2]{N2}). 
For any $n,m\in\N$ with $1\leq n\leq N$, 
there exist $R\in\N$, 
natural numbers $m^{(1)},m^{(2)},\dots,m^{(R)}\geq m$ and 
a family of projections 
\[ e^{(r)}_j \qquad (r=1,2,\dots,R, \ j=0,1,\dots,m^{(r)}-1) \]
in $A_\omega=A^\omega\cap A'$ such that 
\[ \sum_{r=1}^R\sum_{j=0}^{m^{(r)}-1}e^{(r)}_j=1, 
\quad \alpha_{\xi_n}(e^{(r)}_j)=e^{(r)}_{j+1}, 
\quad \alpha_{\xi_i}(e^{(r)}_j)=e^{(r)}_j \]
for any $r=1,2,\dots,R$, $i=1,2,\dots,N$ with $i\neq n$ and 
$j=0,1,\dots,m^{(r)}-1$, 
where the index is understood modulo $m^{(r)}$. 
\end{rem}

\begin{rem}\label{restriction}
It is also obvious that 
if $\alpha$ is an action of $\Z^N$ on $A$ with the Rohlin property, 
then the action $\alpha'$ of $\Z^{N-1}$ generated by 
$\alpha_{\xi_2},\alpha_{\xi_3},\dots,\alpha_{\xi_N}$ also has 
the Rohlin property as a $\Z^{N-1}$-action.  
\end{rem}

\section{Cohomology vanishing}

Throughout this section, 
we let $A$ denote a $C^*$-algebra which is isomorphic to $\mathcal{O}_2$. 
First we need a technical lemma about homotopies of unitaries. 

\begin{lem}\label{HaagerupRordam}
Let $(X,d)$ be a compact metric space and 
let $z:X\to\mathcal{U}(A^\omega)$ be a map. 
Suppose that there exists $C>0$ such that 
$\lVert z(x)-z(x')\rVert\leq Cd(x,x')$ for any $x,x'\in X$. 
Then, for any separable subset $B$ of $A^\omega$, 
we can find a map $\tilde{z}:X\times[0,1]\to\mathcal{U}(A^\omega)$ 
such that the following are satisfied. 
\begin{enumerate}
\item For any $x\in X$, $\tilde{z}(x,0)=z(x)$ and $\tilde{z}(x,1)=1$. 
\item For any $x,x'\in X$ and $t,t'\in[0,1]$, 
\[ \lVert\tilde{z}(x,t)-\tilde{z}(x',t')\rVert
\leq 4Cd(x,x')+\frac{8\pi}{3}\lvert t-t'\rvert. \]
\item For any $a\in B$ and $(x,t)\in X\times[0,1]$, 
$\lVert[\tilde{z}(x,t),a]\rVert\leq4\lVert[z(x),a]\rVert$. 
\end{enumerate}
\end{lem}
\begin{proof}
Since $\mathcal{O}_2$ is isomorphic to 
the infinite tensor product $\bigotimes_{i=1}^\infty\mathcal{O}_2$ 
(see \cite{R2} or \cite[Theorem 3.8]{KP}), 
there exists a unital $C^*$-subalgebra $D$ 
in $A^\omega\cap(B\cup\{z(x)\mid x\in X\})'$ 
such that $D\cong\mathcal{O}_2$. 
We regard $z$ as a unitary in $C(X)\otimes A^\omega$. 
By \cite[Lemma 5.1]{HR} and its proof (see also \cite[Lemma 6]{N2}), 
we can find a unitary $\tilde{z}\in C(X)\otimes C([0,1])\otimes A^\omega$ 
such that the following hold. 
\begin{itemize}
\item $\tilde{z}(x,0)=z(x)$ and $\tilde{z}(x,1)=1$ for all $x\in X$. 
\item $\displaystyle\lVert\tilde{z}(x,t)-\tilde{z}(x,t')\rVert
\leq\frac{8\pi}{3}\lvert t-t'\rvert$ for all $x\in X$ and $t,t'\in[0,1]$. 
\item $\lVert\tilde{z}(x,t)-\tilde{z}(x',t)\rVert
\leq4\lVert z(x)-z(x')\rVert$ for all $x,x'\in X$ and $t\in[0,1]$. 
\item $\lVert[\tilde{z}(x,t),a]\rVert\leq4\lVert[z(x),a]\rVert$ 
for all $(x,t)\in X\times[0,1]$ and $a\in A$. 
\end{itemize}
Then the conclusions follow immediately. 
\end{proof}

Let $N$ be a natural number. 
We denote the $l^\infty$-norm on $\R^N$ by $\lVert\cdot\rVert$. 
We put 
\[ E=\{t\in\R^N\mid\lVert t\rVert\leq1\} \]
and 
\[ \partial E=\{t\in\R^N\mid\lVert t\rVert=1\}. \]

\begin{lem}\label{extendz}
Let $C>0$ and 
let $z_0:\partial E\to\mathcal{U}(A^\omega)$ be a map 
such that 
$\lVert z_0(t)-z_0(t')\rVert\leq C\lVert t-t'\rVert$ 
for every $t,t'\in\partial E$. 
Then, for any separable subset $B$ of $A^\omega$, 
there exists a map $z:E\to\mathcal{U}(A^\omega)$ 
such that the following hold. 
\begin{enumerate}
\item For $t\in\partial E$, $z(t)=z_0(t)$. 
\item For every $t,t'\in E$, 
$\lVert z(t)-z(t')\rVert\leq(24C+16\pi/3)\lVert t-t'\rVert$. 
\item For any $a\in B$ and $t\in E$, 
$\lVert[z(t),a]\rVert\leq4\sup\{\lVert[z(s),a]\rVert\mid s\in\partial E\}$. 
\end{enumerate}
\end{lem}
\begin{proof}
Lemma \ref{HaagerupRordam} applies and 
yields a map $\tilde{z}_0:\partial E\times[0,1]\to\mathcal{U}(A^\omega)$. 
We define $z:E\to\mathcal{U}(A^\omega)$ by 
\[ z(t)=\begin{cases}
1 & \text{ if }\lVert t\rVert\leq1/2 \\
\tilde{z}_0(t/\lVert t\rVert,2(1-\lVert t\rVert)) 
& \text{ if }\lVert t\rVert\geq1/2. \end{cases} \]
Conditions (1) and (3) are immediate from the definition. 

In order to verify (2), take $t,t'\in E$. 
If $\lVert t\rVert,\lVert t'\rVert\leq1/2$, 
we have nothing to do. 
Let us consider the case such that 
$\lVert t'\rVert\leq1/2\leq\lVert t\rVert$. 
Since 
\[ \lVert t\rVert-1/2\leq
\lVert t\rVert-\lVert t'\rVert
\leq\lVert t-t'\rVert, \]
we have 
\begin{align*}
\lVert z(t)-z(t')\rVert
&=\lVert
\tilde{z}_0(t/\lVert t\rVert,2(1-\lVert t\rVert))-1
\rVert \\
&=\lVert
\tilde{z}_0(t/\lVert t\rVert,2(1-\lVert t\rVert))
-\tilde{z}_0(t/\lVert t\rVert,1)
\rVert \\
&\leq\frac{8\pi}{3}
\lvert2(1-\lVert t\rVert)-1\rvert \\
&\leq\frac{16\pi}{3}\lVert t-t'\rVert. 
\end{align*}

It remains for us to check the case such that 
$\lVert t\rVert,\lVert t'\rVert\geq1/2$. 
Put $s=\displaystyle\frac{\lVert t'\rVert}{\lVert t\rVert}t$. 
Then 
\begin{align*}
\lVert z(t)-z(s)\rVert
&\leq\frac{8\pi}{3}
\lvert2(1-\lVert t\rVert)-2(1-\lVert t'\rVert)\rvert \\
&\leq\frac{16\pi}{3}
\lvert\lVert t\rVert-\lVert t'\rVert\rvert \\
&\leq\frac{16\pi}{3}\lVert t-t'\rVert. 
\end{align*}
Besides, 
\begin{align*}
\lVert z(s)-z(t')\rVert
&\leq4C\left\lVert
\frac{t}{\lVert t\rVert}-\frac{t'}{\lVert t'\rVert}
\right\rVert \\
&\leq24C\lVert t-t'\rVert. 
\end{align*}
Combining these, we obtain 
\[ \lVert z(t)-z(t')\rVert
\leq(24C+16\pi/3)\lVert t-t'\rVert. \]
\end{proof}

For each $i=1,2,\dots,N$, we let 
\[ E^+_i=\{(t_1,t_2,\dots,t_N)\in E\mid t_i=1\} \]
and 
\[ E^-_i=\{(t_1,t_2,\dots,t_N)\in E\mid t_i=-1\}. \]
Thus, $E^+_i$ and $E^-_i$ are $(N-1)$-dimensional faces of $E$. 
Let $\sigma_i:E\to E$ be the map such that 
\[ \sigma_i:(t_1,t_2,\dots,t_i,\dots,t_N)
\mapsto(t_1,t_2,\dots,-t_i,\dots,t_N) \]
for each $i=1,2,\dots,N$. 
For $k=1,2,\dots,N-1$, we define $E(k)$ by 
\[ E(k)=\{(t_1,t_2,\dots,t_N)\in E\mid
\#\{i\mid\lvert t_i\rvert=1\}\geq N-k\}. \]
In other words, $E(k)$ is the union of $k$-dimensional faces of $E$. 
We have $\partial E=E(N-1)=\bigcup_{i=1}^N(E^+_i\cup E^-_i)$. 
For each $I\subset\{1,2,\dots,N\}$, we let 
\[ E^+(I)=\{(t_1,t_2,\dots,t_N)\in E\mid
t_i=1\text{ for all }i\notin I\}. \]
Thus, $E^+(I)$ is a $\#I$-dimensional face of $E$ and 
$E^+(I)\cap E(\#I-1)$ is the boundary of $E^+(I)$. 
We also have $E^+_i=E^+(\{1,2,\dots,N\}\setminus\{i\})$. 

\begin{lem}\label{extendz2}
Let $\alpha_1,\alpha_2,\dots,\alpha_N$ be 
$N$ commuting automorphisms on $A^\omega$. 
Suppose that 
there exists a family of unitaries $u_1,u_2,\dots,u_N$ in $A^\omega$ 
satisfying $u_i\alpha_i(u_j)=u_j\alpha_j(u_i)$ for all $i,j=1,2,\dots,N$. 
Let $1\leq k\leq N-2$ and $C>0$. 
Suppose that 
$z_0:E(k)\to\mathcal{U}(A^\omega)$ satisfies the following. 
\begin{itemize}
\item $z_0(1,1,\dots,1)=1$. 
\item For every $i=1,2,\dots,N$ and $t\in E(k)\cap E^+_i$, 
$z_0(\sigma_i(t))=u_i\alpha_i(z_0(t))$. 
\item For every $t,t'\in E(k)$, 
$\lVert z_0(t)-z_0(t')\rVert\leq C\lVert t-t'\rVert$. 
\end{itemize}
Then, for any separable subset $B$ of $A^\omega$, we can find 
an extension $z:E(k+1)\to\mathcal{U}(A^\omega)$ of $z_0$ 
such that the following hold. 
\begin{enumerate}
\item For every $i=1,2,\dots,N$ and $t\in E(k+1)\cap E^+_i$, 
$z(\sigma_i(t))=u_i\alpha_i(z(t))$. 
\item For every $t,t'\in E(k+1)$, 
$\lVert z(t)-z(t')\rVert\leq(48C+32\pi/3)\lVert t-t'\rVert$. 
\item Let $I$ be a subset of $\{1,2,\dots,N\}$ such that $\#I=k+1$. 
For every $a\in B$ and $t\in E^+(I)$, one has 
\[ \lVert[z(t),a]\rVert
\leq4\sup\{\lVert[z_0(s),a]\rVert,
\lVert[z_0(s),\alpha_i^{-1}(a)]\rVert+\lVert[u_i,a]\rVert
\mid s\in E^+(I\setminus\{i\}),i\in I\}. \]
\end{enumerate}
\end{lem}
\begin{proof}
For each $I\subset\{1,2,\dots,N\}$ with $\#I=k+1$, 
by using Lemma \ref{extendz}, 
we can extend $z_0$ on $E^+(I)\cap E(k)$ 
to the map $z_I:E^+(I)\to\mathcal{U}(A^\omega)$ satisfying the following. 
\begin{itemize}
\item For every $t,t'\in E^+(I)$, 
$\lVert z_I(t)-z_I(t')\rVert\leq(24C+16\pi/3)\lVert t-t'\rVert$. 
\item For any $a\in B$ and $t\in E^+(I)$, 
$\lVert[z_I(t),a]\rVert
\leq4\sup\{\lVert z_0(s),a]\rVert\mid s\in E^+(I)\cap E(k)\}$. 
\end{itemize}
We define $z:E(k+1)\to\mathcal{U}(A^\omega)$ as follows. 
First, for $t\in E^+(I)$, we let $z(t)=z_I(t)$. 
Then, we can uniquely extend $z$ to $E(k+1)$ so that 
$z(\sigma_i(t))=u_i\alpha_i(z(t))$ holds 
for any $i=1,2,\dots,N$ and $t\in E(k+1)\cap E^+_i$, 
because of the equality $u_i\alpha_i(u_j)=u_j\alpha_j(u_i)$. 
Note that 
if $t,t'\in E(k+1)$ lie on the same $(k+1)$-dimensional face of $E$, 
then we still have 
$\lVert z(t)-z(t')\rVert\leq(24C+16\pi/3)\lVert t-t'\rVert$. 

Condition (1) is already achieved. 
Let us verify (2). 
Take $t=(t_1,t_2,\dots,t_N)$ and $t'=(t'_1,t'_2,\dots,t'_N)$ in $E(k+1)$. 
Since any unitaries are within distance two of each other, 
we may assume that $\lVert t-t'\rVert$ is less than two. 
We define $s=(s_1,s_2,\dots,s_N)\in E(k+1)$ by 
\[ s_i=\begin{cases}
t_i & \text{ if }\lvert t_i\rvert\neq1\text{ and }\lvert t'_i\rvert\neq1 \\
t_i & \text{ if }\lvert t_i\rvert=1 \\
t'_i & \text{ if }\lvert t'_i\rvert=1. \end{cases} \]
It is easy to see that 
$t$ and $s$ lie on the same $(k+1)$-dimensional face of $E$ 
and that $t'$ and $s$ lie on the same $(k+1)$-dimensional face of $E$. 
In addition, 
both $\lVert t-s\rVert$ and $\lVert s-t'\rVert$ 
are less than $\lVert t-t'\rVert$. 
It follows that 
\begin{align*}
\lVert z(t)-z(t')\rVert
&\leq\lVert z(t)-z(s)\rVert+\lVert z(s)-z(t')\rVert \\
&\leq(24C+16\pi/3)(\lVert t-s\rVert+\lVert s-t'\rVert) \\
&\leq(48C+32\pi/3)\lVert t-t'\rVert, 
\end{align*}
which ensures condition (2). 

Let us consider (3). Take $t\in E^+(I)$. 
We already have 
\[ \lVert[z_I(t),a]\rVert
\leq4\sup\{\lVert z_0(s),a]\rVert\mid s\in E^+(I)\cap E(k)\}. \]
For each $s=(s_1,s_2,\dots,s_N)$ in $E^+(I)\cap E(k)$, 
there exists $i\in I$ such that $\lvert s_i\rvert=1$. 
If $s_i=1$, then $s$ is in $E^+(I\setminus\{i\})$. 
If $s_i=-1$, then $s'=\sigma_i(s)$ is in $E^+(I\setminus\{i\})$ and 
\begin{align*}
\lVert[z_0(s),a]\rVert
&=\lVert[z_0(\sigma_i(s')),a]\rVert \\
&=\lVert[u_i\alpha_i(z_0(s')),a]\rVert \\
&\leq\lVert[z_0(s'),\alpha_i^{-1}(a)]\rVert+\lVert[u_i,a]\rVert, 
\end{align*}
thereby completing the proof. 
\end{proof}

The following proposition is a crucial tool for cohomology vanishing. 

\begin{prop}\label{Berg's}
Let $\alpha_1,\alpha_2,\dots,\alpha_N\in\Aut(A^\omega)$ be 
$N$ commuting automorphisms on $A^\omega$. 
Suppose that 
there exists a family of unitaries $u_1,u_2,\dots,u_N$ in $A^\omega$ 
satisfying $u_i\alpha_i(u_j)=u_j\alpha_j(u_i)$ for all $i,j=1,2,\dots,N$. 
Then, for any separable subset $B$ of $A^\omega$, 
we can find a continuous map $z:E\to\mathcal{U}(A^\omega)$ 
such that the following hold. 
\begin{enumerate}
\item $z(1,1,\dots,1)=1$. 
\item For every $i=1,2,\dots,N$ and $t\in E^+_i$, 
$z(\sigma_i(t))=u_i\alpha_i(z(t))$. 
\item For every $t,t'\in E$, 
$\lVert z(t)-z(t')\rVert\leq50^N\lVert t-t'\rVert$. 
\item For any $a\in B$ and $t\in E$, 
\[ \lVert[z(t),a]\rVert
\leq4^N\sup\sum_{k=1}^K
\lVert[u_{i_k},
(\alpha_{i_1}\alpha_{i_2}\dots\alpha_{i_{k-1}})^{-1}(a)]\rVert, \]
where the supremum is taken over 
all permutations $i_1,i_2,\dots,i_K$ of elements in $\{1,2,\dots,N\}$. 
\end{enumerate}
\end{prop}
\begin{proof}
Clearly we may assume that 
$B$ is $\alpha_i$-invariant for every $i=1,2,\dots,N$. 
By applying Lemma \ref{HaagerupRordam} to the case 
such that $X$ is a singleton, 
for each $i=1,2,\dots,N$, we obtain 
a map $z_{1,i}$ from $E^+(\{i\})\cong[0,1]$ to $\mathcal{U}(A^\omega)$ 
satisfying the following. 
\begin{itemize}
\item $z_{1,i}(\sigma_i(1,1,\dots,1))=u_i$ and $z_{1,i}(1,1,\dots,1)=1$. 
\item For any $t,t'\in E^+(\{i\})$, 
$\displaystyle \lVert z_{1,i}(t)-z_{1,i}(t')\rVert
\leq\frac{8\pi}{3}\lVert t-t'\rVert$. 
\item For any $a\in B$ and $t\in E^+(\{i\})$, 
$\lVert[z_{1,i}(t),a]\rVert\leq4\lVert[u_i,a]\rVert$. 
\end{itemize}
From these maps $z_{1,i}$'s, 
by the same argument as in the previous lemma, 
we can construct a map $z_1:E(1)\to\mathcal{U}(A^\omega)$ 
such that the following are satisfied. 
\begin{itemize}
\item For every $i=1,2,\dots,N$ and $t\in E^+(\{i\})$, 
$z_1(t)=z_{1,i}(t)$. 
\item For every $i=1,2,\dots,N$ and $t\in E(1)\cap E^+_i$, 
$z_1(\sigma_i(t))=u_i\alpha_i(z_1(t))$. 
\item For every $t,t'\in E(1)$, 
$\displaystyle \lVert z_1(t)-z_1(t')\rVert
\leq\frac{16\pi}{3}\lVert t-t'\rVert$. 
\end{itemize}
Note that we have used the equality $u_i\alpha_i(u_j)=u_j\alpha_j(u_i)$. 

We apply Lemma \ref{extendz2} to $z_1:E(1)\to\mathcal{U}(A^\omega)$ 
and obtain an extension $z_2:E(2)\to\mathcal{U}(A^\omega)$ of $z_1$ 
which satisfies the following. 
\begin{itemize}
\item For every $i=1,2,\dots,N$ and $t\in E(2)\cap E^+_i$, 
$z_2(\sigma_i(t))=u_i\alpha_i(z_2(t))$. 
\item For every $t,t'\in E(2)$, 
$\lVert z_2(t)-z_2(t')\rVert\leq50^2\lVert t-t'\rVert$. 
\item Let $I$ be a subset of $\{1,2,\dots,N\}$ such that $\#I=2$. 
For every $a\in B$ and $t\in E^+(I)$, one has 
\[ \lVert[z_2(t),a]\rVert
\leq4\sup\{\lVert[z_1(s),a]\rVert,
\lVert[z_1(s),\alpha_i^{-1}(a)]\rVert+\lVert[u_i,a]\rVert
\mid s\in E^+(I\setminus\{i\}),i\in I\}. \]
\end{itemize}
Repeating this argument, we get a map $z_{N-1}$ 
from $E(N-1)=\partial E$ to $\mathcal{U}(A^\omega)$ 
satisfying the following. 
\begin{itemize}
\item For every $i=1,2,\dots,N$ and $t\in E^+_i$, 
$z_{N-1}(\sigma_i(t))=u_i\alpha_i(z_{N-1}(t))$. 
\item For every $t,t'\in E(N-1)$, 
$\lVert z_{N-1}(t)-z_{N-1}(t')\rVert\leq50^{N-1}\lVert t-t'\rVert$. 
\item Let $I$ be a subset of $\{1,2,\dots,N\}$ such that $\#I=N-1$. 
For every $a\in B$ and $t\in E^+(I)$, one has 
\begin{align*}
&\lVert[z_{N-1}(t),a]\rVert \\
&\leq4\sup\{\lVert[z_{N-2}(s),a]\rVert,
\lVert[z_{N-2}(s),\alpha_i^{-1}(a)]\rVert+\lVert[u_i,a]\rVert
\mid s\in E^+(I\setminus\{i\}),i\in I\}. 
\end{align*}
\end{itemize}
By using Lemma \ref{extendz}, 
we can extend $z_{N-1}$ to $z:E\to\mathcal{U}(A^\omega)$. 
Then, clearly conditions (2) and (3) are satisfied. 
As for condition (4), we have 
\begin{align*}
\lVert[z(t),a]\rVert
&\leq4\sup\{\lVert[z_{N-1}(s),a]\rVert\mid s\in E(N-1)\} \\
&\leq4\sup\{\lVert[z_{N-1}(s),a]\rVert,
\lVert[z_{N-1}(s),\alpha_i^{-1}(a)]\rVert+\lVert[u_i,a]\rVert
\mid s\in E^+_i,1\leq i\leq N\}
\end{align*}
for any $t\in E$ and $a\in B$. 
By estimating norms of 
commutators of $z_k(s)$ with elements in $B$ inductively, 
we obtain the desired inequality. 
\end{proof}

\bigskip

Now we would like to show the cohomology vanishing theorem. 

\begin{thm}\label{CVanish}
Let $\alpha$ be an action of $\Z^N$ on $A$ with the Rohlin property 
and let $\tilde\alpha$ be a perturbed action of $\alpha$ 
on $A^\omega$ by an $\alpha$-cocycle in $A^\omega$. 
Let $B\subset A^\omega$ be an $\tilde\alpha$-invariant separable subset. 
Suppose that a family of unitaries 
$\{u_n\}_{n\in\Z^N}$ in $A^\omega\cap B'$ is an $\tilde\alpha$-cocycle. 
Then, there exists a unitary $v\in\mathcal{U}(A^\omega\cap B')$ 
such that 
\[ u_n=v\tilde\alpha_n(v^*) \]
for each $n\in\Z^N$, that is, $\{u_n\}_{n\in\Z^N}$ is a coboundary. 
\end{thm}
\begin{proof}
Evidently it suffices to show the following: 
For any $\ep>0$, 
there exists a unitary $v\in\mathcal{U}(A^\omega\cap B')$ 
such that 
\[ \lVert u_{\xi_i}-v\tilde\alpha_{\xi_i}(v^*)\rVert<\ep \]
for every $i=1,2,\dots,N$. 
Choose a natural number $M$ so that $\ep(M-1)>2\cdot50^N$.  
Since $\alpha$ has the Rohlin property, 
there exist $R\in\N$ and $m^{(1)},m^{(2)},\dots,m^{(R)}\in\Z^N$ 
with $m^{(1)},\dots,m^{(R)}\geq(M,M,\dots,M)$ 
and which satisfies the requirement in Definition \ref{Rohlin}. 

For each $r=1,2,\dots,R$ and $i=1,2,\dots,N$, 
let $m^{(r)}_i$ be the $i$-th summand of $m^{(r)}$. 
We put $\eta_{r,i}=m^{(r)}_i\xi_i\in\Z^N$. 
By applying Lemma \ref{Berg's} to 
$\tilde\alpha_{\eta_{r,1}},\tilde\alpha_{\eta_{r,2}},
\dots,\tilde\alpha_{\eta_{r,N}}$ and 
unitaries $u_{\eta_{r,1}},u_{\eta_{r,2}},\dots,u_{\eta_{r,N}}$ 
in $\mathcal{U}(A^\omega\cap B')$, 
we obtain a map $z^{(r)}:E\to\mathcal{U}(A^\omega\cap B')$ 
satisfying the following. 
\begin{itemize}
\item $z^{(r)}(1,1,\dots,1)=1$. 
\item For every $i=1,2,\dots,N$ and $t\in E^+_i$, 
$z^{(r)}(\sigma_i(t))
=u_{\eta_{r,i}}\tilde\alpha_{\eta_{r,i}}(z^{(r)}(t))$. 
\item For every $t,t'\in E$, 
$\lVert z^{(r)}(t)-z^{(r)}(t')\rVert\leq50^N\lVert t-t'\rVert$. 
\end{itemize}
For each $r=1,2,\dots,R$ and $g=(g_1,g_2,\dots,g_N)\in\Z_{m^{(r)}}$, 
we define $w^{(r)}_g$ in $\mathcal{U}(A^\omega\cap B')$ by 
\[ w^{(r)}_g=z^{(r)}\left(\frac{2g_1}{m^{(r)}_1-1}-1,
\frac{2g_2}{m^{(r)}_2-1}-1,\dots,\frac{2g_N}{m^{(r)}_N-1}-1\right). \]
It is easily seen that 
one has the following 
for any $r=1,2,\dots,R$, $i=1,2,\dots,N$ 
and $g=(g_1,g_2,\dots,g_N)\in\Z_{m^{(r)}}$. 
\begin{itemize}
\item If $g_i\neq0$, then 
$\lVert w^{(r)}_g-w^{(r)}_{g-\xi_i}\rVert$ is less than $\ep$. 
\item If $g_i=0$, then 
$w^{(r)}_g$ is equal to 
$u_{\eta_{r,i}}\tilde\alpha_{\eta_{r,i}}(w^{(r)}_{g+\eta_{r,i}-\xi_i})$. 
\end{itemize}
By Remark \ref{reindex}, 
we can take a family of projections 
$\{e^{(r)}_g\mid r=1,2,\dots,R, \ g\in\Z_{m^{(r)}}\}$ 
in $A^\omega\cap B'$ such that 
\[ \sum_{r=1}^R\sum_{g\in\Z_{m^{(r)}}}e^{(r)}_g=1, 
\quad \tilde\alpha_{\xi_i}(e^{(r)}_g)=e^{(r)}_{g+\xi_i} \]
for any $r=1,2,\dots,R$, $i=1,2,\dots,N$ and $g\in\Z_{m^{(r)}}$, 
where $g+\xi_i$ is understood modulo $m^{(r)}\Z^N$. 
Moreover, we may also assume that 
$e^{(r)}_g$ commutes with $u_g$ and $\tilde\alpha_g(w^{(r)}_g)$. 
Define $v\in\mathcal{U}(A^\omega\cap B')$ by 
\[ v=\sum_{r=1}^R\sum_{g\in\Z_{m^{(r)}}}
u_g\tilde\alpha_g(w^{(r)}_g)e^{(r)}_g. \]
It is now routinely checked that 
$\lVert u_{\xi_i}-v\tilde\alpha_{\xi_i}(v^*)\rVert$ is 
less than $\ep$ for each $i=1,2,\dots,N$. 
\end{proof}

The following corollary is an immediate consequence of the theorem above 
and we omit the proof. 

\begin{cor}\label{appCVanish}
Let $\alpha$ be an action of $\Z^N$ on $A$ with the Rohlin property. 
For any $\ep>0$ and a finite subset $\mathcal{F}$ of $A$, 
there exist $\delta>0$ and a finite subset $\mathcal{G}$ of $A$ 
such that the following holds: 
If a family of unitaries $\{u_n\}_{n\in\Z^N}$ in $A$ 
is an $\alpha$-cocycle satisfying 
\[ \lVert[u_{\xi_i},a]\rVert<\delta \]
for every $i=1,2,\dots,N$ and $a\in\mathcal{G}$, then 
we can find a unitary $v\in\mathcal{U}(A)$ 
satisfying 
\[ \lVert u_{\xi_i}-v\alpha_{\xi_i}(v^*)\rVert<\ep \]
and 
\[ \lVert[v,a]\rVert<\ep \]
for each $i=1,2,\dots,N$ and $a\in\mathcal{F}$. 
Furthermore, 
if $\mathcal{F}$ is an empty set, 
then we can take an empty set for $\mathcal{G}$. 
\end{cor}

\section{Rohlin type theorem}

Throughout this section, 
we let $A$ denote a unital $C^*$-algebra 
which is isomorphic to $\mathcal{O}_2$. 
In this section, we would like to show the Rohlin type theorem 
for $\Z^N$-actions on $A$ 
by combining techniques developed in \cite{N1} and \cite{N2}. 

\begin{lem}\label{scattered}
Let $\alpha$ be an outer action of $\Z^N$ on $A$. 
Then, for any $m\in\N^N$ and 
a non-zero projection $p\in A_\omega=A^\omega\cap A'$, 
there exists a non-zero projection $e$ in $pA_\omega p$ such that 
$e\alpha_g(e)=0$ for all $g\in\Z_m\setminus\{0\}$. 
\end{lem}
\begin{proof}
We can prove this lemma 
exactly in the same way as \cite[Lemma 3]{N2}. 
See also \cite[Lemma 2]{N2}. 
\end{proof}

The following lemma is a generalization of \cite[Lemma 3.5]{K2}. 

\begin{lem}\label{fixedisom}
Let $\alpha$ be an action of $\Z^N$ on $A$ with the Rohlin property. 
Suppose that 
one has two non-zero projections $e,f\in A_\omega$ 
satisfying $\alpha_n(e)=e$ and $\alpha_n(f)=f$ for any $n\in\Z^N$. 
Then, there exists $w\in A_\omega$ such that 
$w^*w=e$, $ww^*=f$ and $\alpha_n(w)=w$ for all $n\in\Z^N$. 
\end{lem}
\begin{proof}
To simplify notation, we denote $\alpha_{\xi_i}$ by $\alpha_i$. 
Since $A_\omega$ is purely infinite simple and $K_0(A_\omega)=0$, 
there exists a partial isometry $u\in A_\omega$ 
such that $u^*u=e$, $uu^*=f$. 
Put $u_i=u^*\alpha_i(u)+1-e$. 
It is straightforward to verify $u_i\alpha_i(u_j)=u_j\alpha_j(u_i)$ 
for all $i,j=1,2,\dots,N$. 
Thus, the family $\{u_i\}$ is an $\alpha$-cocycle in $A_\omega$. 
Clearly $u_i$ commutes with $e$. 
By Theorem \ref{CVanish}, 
we can find a unitary $v\in A_\omega$ such that 
$[v,e]=0$ and $u_i=v\alpha_i(v^*)$. 
Then $w=uv$ satisfies the requirements. 
\end{proof}

We also have an approximate version of the lemma above. 

\begin{lem}\label{appfixedisom}
Let $\alpha$ be an action of $\Z^N$ on $A$ with the Rohlin property. 
For any $\ep>0$ and a finite subset $\mathcal{F}$ of $A$, 
there exist $\delta>0$ and a finite subset $\mathcal{G}$ of $A$ 
such that the following holds: 
Suppose that 
$e$ and $f$ are two non-zero projections in $A$ satisfying 
\[ \lVert[e,a]\rVert<\delta, \ \lVert[f,a]\rVert<\delta
\quad\text{ for all }a\in\mathcal{G} \]
and 
\[ \lVert\alpha_{\xi_i}(e)-e\rVert<\delta, \ 
\lVert\alpha_{\xi_i}(f)-f\rVert<\delta
\quad\text{ for each }i=1,2,\dots,N. \]
Then, we can find a partial isometry $v\in A$ such that 
$v^*v=e$, $vv^*=f$ and 
\[ \lVert[v,a]\rVert<\ep
\quad\text{ for all }a\in\mathcal{F} \]
and 
\[ \lVert\alpha_{\xi_i}(v)-v\rVert<\ep
\quad\text{ for each }i=1,2,\dots,N. \]
\end{lem}
\begin{proof}
This immediately follows from the lemma above. 
\end{proof}

Next, we have to recall a technical result 
about almost cyclic projections. 

Suppose that we are given $\ep>0$ and $n\in\N$. 
Choose a natural number $k\in\N$ so that $2/\sqrt{k}<\ep$. 
Let $\alpha$ be an automorphism on $A$ and 
let $p$ be a projection of $A$ 
which satisfies $p\alpha^i(p)=0$ for all $i=1,2,\dots,nk$. 
Furthermore, let $v\in A$ be a partial isometry 
such that $v^*v=p$ and $vv^*=\alpha(p)$. 
Define 
\[ E_{i,j}=\begin{cases}
\alpha^{i-1}(v)\alpha^{i-1}(v)\dots\alpha^j(v) & \text{ if }i>j \\
\alpha^i(p) & \text{ if }i=j \\
\alpha^i(v^*)\alpha^{i+1}(v^*)\dots\alpha^{j-1}(v^*) & 
\text{ if }i<j \end{cases} \]
for each $i,j=0,1,\dots,nk$. 
Then we can easily see that $\{E_{i,j}\}$ is a system of matrix units 
and $\alpha(E_{i,j})=E_{i+1,j+1}$ for any $i,j=0,1,\dots,nk-1$. 
We let 
\[ f=\frac{1}{k}\sum_{i,j=0}^{k-1}E_{ni,nj} \]
and 
\[ e_i=\alpha^i(f) \]
for all $i=0,1,\dots,n-1$. 
Then $\{e_i\}$ is an orthogonal family of projections in $A$ satisfying 
\[ e_0+e_1+\dots+e_{n-1}\leq\sum_{i=0}^{nk-1}E_{i,i} \]
and 
\[ \lVert e_0-\alpha(e_{n-1})\rVert<\frac{2}{\sqrt{k}}<\ep. \]
Moreover, this argument applies to the case that 
$p$ is almost orthogonal to $\alpha^i(p)$. 
Thus, for any $\ep>0$, 
there exists a small positive constant $c(\ep,n,k)$ 
such that the following holds: 
if $\lVert p\alpha^i(p)\rVert<c(\ep,n,k)$ for all $i=1,2,\dots,nk$, 
then we can find a projection $e_0$ such that 
\[ \left\lVert(e_0+\alpha(e_0)+\dots+\alpha^{n-1}(e_0))
\left(1-\sum_{i=0}^{nk-1}E_{i,i}\right)\right\rVert<\ep \]
and 
\[ \lVert e_0-\alpha^n(e_0)\rVert<\frac{2}{\sqrt{k}}<\ep. \]
Using these estimates, we can prove the next lemma. 

\begin{lem}\label{appcyclicproj}
Let $\alpha$ be an outer action of $\Z^N$ on $A$. 
Suppose that 
the action $\alpha'$ of $\Z^{N-1}$ generated 
by $\alpha_{\xi_2},\alpha_{\xi_3},\dots,\alpha_{\xi_N}$ has 
the Rohlin property. 
Then, for any $n\in\N$, $\ep>0$ and 
a finite subset $\mathcal{F}$ of $A$, 
there exist $m\in\N^N$, $\delta>0$ and 
a finite subset $\mathcal{G}$ of $A$ satisfying the following: 
If a non-zero projection $p$ in $A$ satisfies 
\[ \lVert p\alpha_g(p)\rVert<\delta
\quad\text{ for each }g\in\Z_m\setminus\{0\} \]
and 
\[ \lVert[\alpha_g(p),a]\rVert<\delta
\quad\text{ for all }g\in\Z_m\text{ and }a\in\mathcal{G}, \]
then there exists a non-zero projection $e$ in $A$ such that 
\begin{enumerate}
\item $\left\lVert(e+\alpha_{\xi_1}(e)+\dots+\alpha^{n-1}_{\xi_1}(e))
\left(1-\sum_{g\in\Z_m}\alpha_g(p)\right)\right\rVert<\ep$. 
\item $\lVert[\alpha^j_{\xi_1}(e),a]\rVert<\ep$ 
for all $j=0,1,\dots,n-1$ and $a\in\mathcal{F}$. 
\item $\lVert\alpha_{\xi_i}(e)-e\rVert<\ep$ for each $i=2,3,\dots,N$. 
\item $\lVert\alpha^n_{\xi_1}(e)-e\rVert<\ep$. 
\end{enumerate}
\end{lem}
\begin{proof}
We prove the lemma by induction on $N$. 
When $N=1$, the assertion follows immediately 
from \cite[Lemma 4]{N2} and its proof. 
Suppose that the case of $N-1$ has been shown. 
We would like to consider the case of $\Z^N$-actions. 
To simplify notation, we denote $\alpha_{\xi_i}$ by $\alpha_i$. 

We are given $n\in\N$, $\ep>0$ 
and a finite subset $\mathcal{F}$ of $A$. 
We choose $k\in\N$ so that $2/\sqrt{k}<\ep$. 
We will eventually find 
a projection $q$ and a partial isometry $v\in A$ such that 
$v^*v=q$, $vv^*=\alpha_1(q)$ and 
$\lVert q\alpha^j_1(q)\rVert<c(\ep/2,n,k)$ for all $j=0,1,\dots,nk$. 
Then, by the above mentioned technique, 
we will construct a projection $e$ satisfying 
\[ \left\lVert(e+\alpha_1(e)+\dots+\alpha_1^{n-1}(e))
\left(1-\sum_{j=0}^{nk-1}\alpha_1^j(q)\right)\right\rVert<\ep/2 \]
and 
\[ \lVert e-\alpha_1^n(e)\rVert<\frac{2}{\sqrt{k}}<\ep. \]
In this construction, 
we can find $\ep'>0$ and a finite subset $\mathcal{F}'$ of $A$ 
such that the following hold: 
If the projection $q$ and the partial isometry $v$ satisfy 
\[ \lVert[q,a]\rVert<\ep', \ \lVert[v,a]\rVert<\ep'
\quad\text{ for all }a\in\mathcal{F}' \]
and 
\[ \lVert\alpha_i(q)-q\rVert<\ep', \ \lVert\alpha_i(v)-v\rVert<\ep'
\quad\text{ for each }i=2,3,\dots,N, \]
then the obtained projection $e$ satisfies conditions (2) and (3). 

By applying Lemma \ref{appfixedisom} to 
the action $\alpha'$ of $\Z^{N-1}$, $\ep'>0$ and $\mathcal{F}'$, 
we get $\ep''>0$ and a finite subset $\mathcal{F}''$ of $A$. 
We may assume that $\ep''$ is smaller than 
\[ \min\left\{\ep',\frac{c(\ep/2,n,k)}{3},\frac{\ep}{4nk}\right\} \]
and that $\mathcal{F}''$ contains $\mathcal{F}'$. 
By using the induction hypothesis for 
the action $\alpha'$ of $\Z^{N-1}$, $n=1$, 
$\ep''>0$ and $\mathcal{F}''\cup\alpha^{-1}_1(\mathcal{F}'')$ 
(see Remark \ref{restriction}), 
we have $m'\in\N^{N-1}$, $\delta>0$ and a finite subset $\mathcal{G}$. 
We define $m\in\N^N$ by $m=(nk,m')$. 

In order to show that these items meet the requirements, 
let $p$ be a projection in $A$ such that 
\[ \lVert p\alpha_g(p)\rVert<\delta
\quad\text{ for each }g\in\Z_m\setminus\{0\} \]
and 
\[ \lVert[\alpha_g(p),a]\rVert<\delta
\quad\text{ for all }g\in\Z_m\text{ and }a\in\mathcal{G}. \]
Then there exists a projection $q$ in $A$ 
which satisfies the following. 
\begin{itemize}
\item $\left\lVert q(1-\sum_g\alpha_g(p))\right\rVert<\ep''$, 
where the summation runs over all $g=(g_1,g_2,\dots,g_N)\in\Z_m$ 
with $g_1=0$. 
\item $\lVert[p,a]\rVert<\ep''$ and $\lVert[\alpha_1(p),a]\rVert<\ep''$ 
for all $a\in\mathcal{F}''$. 
\item $\lVert\alpha_i(p)-p\rVert<\ep''$ for each $i=2,3,\dots,N$. 
\end{itemize}
In addition, by taking $\delta$ sufficiently small, 
we may assume that the first condition above implies 
$\lVert q\alpha_1^j(q)\rVert<c(\ep/2,n,k)$ for all $j=1,2,\dots,nk$. 
From Lemma \ref{appfixedisom}, 
there exists a partial isometry $v$ such that 
$v^*v=q$, $vv^*=\alpha_1(q)$ and 
\[ \lVert[v,a]\rVert<\ep'
\quad\text{ for all }a\in\mathcal{F}' \]
and 
\[ \lVert\alpha_i(v)-v\rVert<\ep'
\quad\text{ for each }i=2,3,\dots,N. \]
Then, by the choice of $\ep'>0$ and $\mathcal{F}'$, 
the desired projection $e$ is obtained. 
\end{proof}

By using the lemma above, we can show the following. 

\begin{lem}\label{cyclicproj}
Let $\alpha$ be an outer action of $\Z^N$ on $A$. 
Suppose that 
the action $\alpha'$ of $\Z^{N-1}$ generated 
by $\alpha_{\xi_2},\alpha_{\xi_3},\dots,\alpha_{\xi_N}$ has 
the Rohlin property. 
Then, for any $n\in\N$, 
there exist non-zero mutually orthogonal projections 
$e_0,e_1,\dots,e_{n-1}$ in $A_\omega=A^\omega\cap A'$ 
such that 
\[ \alpha_{\xi_i}(e_j)=e_j \]
for each $i=2,3,\dots,N$ and $j=0,1,\dots,n-1$, and 
\[ \alpha_{\xi_1}(e_j)=e_{j+1} \]
for each $j=0,1,\dots,n-1$, where the addition is understood modulo $n$. 
\end{lem}
\begin{proof}
It suffices to show the following. 
For any $n\in\N$, $\ep>0$ and a finite subset $\mathcal{F}$ of $A$, 
there exist non-zero almost mutually orthogonal projections 
$e_0,e_1,\dots,e_{n-1}$ in $A$ such that the following are satisfied. 
\begin{itemize}
\item $\lVert[e_j,a]\rVert<\ep$ 
for any $j=0,1,\dots,n-1$ and $a\in\mathcal{F}$. 
\item $\lVert\alpha_{\xi_i}(e_j)-e_j\rVert<\ep$ 
for each $i=2,3,\dots,N$ and $j=0,1,\dots,n-1$. 
\item $\lVert\alpha_{\xi_1}(e_j)-e_{j+1}\rVert<\ep$ 
for each $j=0,1,\dots,n-1$, where the addition is understood modulo $n$. 
\end{itemize}
Lemma \ref{appcyclicproj} applies to 
$n\in\N$, $\ep>0$ and $\mathcal{F}\subset A$ 
and yields $m\in\N^N$, $\delta>0$ and 
a finite subset $\mathcal{G}$ of $A$. 
From Lemma \ref{scattered}, 
there exists a non-zero projection $p$ in $A$ such that 
$\lVert p\alpha_g(p)\rVert<\delta$ for all $g\in\Z_m\setminus\{0\}$ and 
$\lVert[\alpha_g(p),a]\rVert<\delta$ 
for all $g\in\Z_m$ and $a\in\mathcal{G}$. 
Then, by Lemma \ref{appcyclicproj}, we can find desired projections. 
\end{proof}

The following lemma is an easy exercise and we omit the proof. 

\begin{lem}\label{matrix}
For any $M\in\N$ and $\ep>0$, 
there exists a natural number $L\in\N$ 
such that the following holds. 
If $u,u'$ are two unitaries in the matrix algebra $M_{LM+1}(\C)$ 
satisfying 
\[ \Sp(u)=\{1\}
\cup\left\{\exp\frac{2\pi\sqrt{-1}j}{LM}\mid j=0,1,\dots,LM-1\right\} \]
and 
\begin{align*}
\Sp(u')=&\left\{\exp\frac{2\pi\sqrt{-1}j}{(L-1)M}
\mid j=0,1,\dots,(L-1)M-1\right\} \\
&\cup\left\{\exp\frac{2\pi\sqrt{-1}j}{M+1}\mid j=0,1,\dots,M\right\}
\end{align*}
with multiplicity, 
then there exists a unitary $w$ in $M_{LM+1}(\C)$ such that 
\[ \lVert wuw^*-u'\rVert<\ep. \]
\end{lem}

Now we can prove the Rohlin type theorem 
for $\Z^N$-actions on the Cuntz algebra $\mathcal{O}_2$. 

\begin{thm}\label{Rohlintype}
Let $A$ be a unital $C^*$-algebra 
which is isomorphic to $\mathcal{O}_2$ and 
let $\alpha$ be an action of $\Z^N$ on $A$. 
Then the following are equivalent. 
\begin{enumerate}
\item $\alpha$ is an outer action. 
\item $\alpha$ has the Rohlin property. 
\end{enumerate}
\end{thm}
\begin{proof}
It is obvious that (2) implies (1). 
The implication from (1) to (2) is shown by induction on $N$. 
When $N=1$, the conclusion follows from \cite[Theorem 3.1]{K2}. 
We assume that the assertion has been shown for $N-1$. 
Let us consider the case of $\Z^N$-actions. 
To simplify notation, we denote $\alpha_{\xi_i}$ by $\alpha_i$. 

It suffices to show the following (see Remark \ref{Rem2ofN1}): 
For any $M\in\N$ and $\ep>0$, there exist projections 
$p_0,p_1,\dots,p_{M-1}$ and 
$q_0,q_1,\dots,q_M$ in $A_\omega$ such that the following hold. 
\begin{itemize}
\item $\sum_{j=0}^{M-1}p_j+\sum_{j=0}^Mq_j=1$. 
\item $\alpha_i(p_j)=p_j$ and $\alpha_i(q_j)=q_j$ 
for each $i=2,3,\dots,N$ and $j=0,1,\dots,M-1$. 
\item $\lVert\alpha_1(p_j)-p_{j+1}\rVert<\ep$ 
for each $j=0,1,\dots,M-1$, 
where $p_M$ is understood as $p_0$. 
\item $\lVert\alpha_1(q_j)-q_{j+1}\rVert<\ep$ 
for each $j=0,1,\dots,M$, 
where $q_{M+1}$ is understood as $q_0$. 
\end{itemize}
Let $\alpha'$ be the action of $\Z^{N-1}$ 
generated by $\alpha_2,\alpha_3,\dots,\alpha_N$. 
By the induction hypothesis, $\alpha'$ has the Rohlin property. 
We are given $M\in\N$ and $\ep>0$. 
By applying Lemma \ref{matrix}, we obtain a natural number $L\in\N$. 
Let $K$ be a very large natural number. 
By Lemma \ref{cyclicproj}, 
we can find non-zero projections $e_0,e_1,\dots,e_{LM-1}$ in $A_\omega$ 
such that 
\[ \alpha_i(e_j)=e_j \]
for each $i=2,3,\dots,N$ and $j=0,1,\dots,LM-1$, and 
\[ \alpha_1(e_j)=e_{j+1} \]
for each $j=0,1,\dots,LM-1$, where the addition is understood modulo $LM$. 
Set $e=\sum_{j=0}^{LM-1}e_j$. 
If $e=1$, we have nothing to do, so we assume $e\neq1$. 

Take $\ep_0>0$ and a finite subset $\mathcal{F}_0$ of $A$ arbitrarily. 
By using Lemma \ref{appcyclicproj} 
for $\alpha'$, $n=1$, $\ep_0$ and $\mathcal{F}_0$, 
we get $m'\in\N^{N-1}$, $\delta>0$ and a finite subset $\mathcal{G}$. 
Define $m\in\N^N$ by $m=(KLM,m')$. 
By Lemma \ref{scattered}, 
we can find a non-zero projection $p$ in $e_0A_\omega e_0$ 
such that $p\alpha_g(p)=0$ for all $g\in\Z_m\setminus\{0\}$. 
Then, by applying Lemma \ref{appcyclicproj} 
to each coordinate of 
a representing sequence of $p$ in $\ell^\infty(\N,A)$, 
we can construct a non-zero projection $f$ in $A^\omega$ 
satisfying the following. 
\begin{itemize}
\item $\lVert f(1-\sum_g\alpha_g(p))\rVert<\ep_0$, 
where the summation runs over all $g=(g_1,g_2,\dots,g_N)\in\Z_m$ 
with $g_1=0$. 
\item $\lVert[f,a]\rVert<\ep_0$ for all $a\in\mathcal{F}_0$. 
\item $\lVert\alpha_i(f)-f\rVert<\ep_0$ for each $i=2,3,\dots,N$. 
\end{itemize}
Notice that, from the first condition, 
one has $\lVert f(1-e_0)\rVert<\ep_0$ 
and $\lVert f\alpha_1^j(f)\rVert<2\ep_0$ for all $j=1,2,\dots,KLM-1$. 
Since $\ep_0>0$ and $\mathcal{F}_0$ were arbitrary, 
by the reindexation trick, we may assume that 
there exists a non-zero projection $f\in A_\omega$ such that 
$f\leq e_0$, 
\[ \alpha_i(f)=f\quad\text{ for each }i=2,3,\dots,N \]
and 
\[ f\alpha_1^j(f)=0\quad\text{ for each }j=1,2,\dots,KLM-1. \]
By applying Lemma \ref{fixedisom} to $\alpha'$, $1-e$ and $f$ 
(note that $\alpha'$ has the Rohlin property), 
we obtain a partial isometry $v\in A_\omega$ such that 
$v^*v=1-e$, $vv^*=f$ and $\alpha_i(v)=v$ for all $i=2,3,\dots,N$. 
We let 
\[ w=\frac{1}{\sqrt{K}}\sum_{k=0}^{K-1}\alpha_1^{kLM}(v). \]
Then, $w$ is a partial isometry in $A_\omega$ 
satisfying $w^*w=1-e$ and $ww^*\leq e_0$. 
In addition, $\alpha_i(w)=w$ for all $i=2,3,\dots,N$ 
and $\lVert\alpha_1^{LM}(w)-w\rVert<2/\sqrt{K}$. 

We consider a $C^*$-subalgebra $D$ of $A_\omega$, 
defined by 
\[ D=C^*(w,\alpha_1(w),\dots,\alpha_1^{LM-1}(w)). \]
It is easy to see that 
$D$ is isomorphic to the matrix algebra $M_{LM+1}(\C)$ and 
its unit is 
\[ 1_D=1-e+ww^*+\alpha_1(ww^*)+\dots+\alpha_1^{LM-1}(ww^*). \]
By choosing $K$ so large, we may assume that 
there exists a unitary 
\[ u=\begin{bmatrix}
1 & & & & & & \\
 & 0 & \cdot & \ldots & \cdot & 0 & 1 \\
 & 1 & 0 & \ldots & \cdot & & 0 \\
 & 0 & 1 & \cdot & & & \cdot \\
 & \vdots & & \ddots & \ddots & \vdots \\
 & \cdot & & 0 & 1 & 0 & 0 \\
 & 0 & & & 0 & 1 & 0
\end{bmatrix} \]
in $D$ such that 
\[ \lVert\alpha_1(x)-uxu^*\rVert\leq\ep\lVert x\rVert \]
for all $x\in D$. 
By the choice of $L$, 
we can find a unitary $u'$ in $D$ such that 
\begin{align*}
\Sp(u')=&\left\{\exp\frac{2\pi\sqrt{-1}j}{(L-1)M}
\mid j=0,1,\dots,(L-1)M-1\right\} \\
&\cup\left\{\exp\frac{2\pi\sqrt{-1}j}{M+1}\mid j=0,1,\dots,M\right\}
\end{align*}
and $\lVert u-u'\rVert<\ep$. 
It follows that 
there exist non-zero mutually orthogonal projections 
$p_0,p_1,\dots,p_{M-1},q_0,q_1,\dots,q_M$ in $D$ satisfying the following. 
\begin{itemize}
\item $\sum_{j=0}^{M-1}p_j+\sum_{j=0}^Mq_j=1_D$. 
\item $\alpha_i(p_j)=p_j$ and $\alpha_i(q_j)=q_j$ 
for each $i=2,3,\dots,N$ and $j=0,1,\dots,M-1$. 
\item $\lVert\alpha_1(p_j)-p_{j+1}\rVert<3\ep$ 
for each $j=0,1,\dots,M-1$, 
where $p_M$ is understood as $p_0$. 
\item $\lVert\alpha_1(q_j)-q_{j+1}\rVert<3\ep$ 
for each $j=0,1,\dots,M$, 
where $q_{M+1}$ is understood as $q_0$. 
\end{itemize}
Finally, we define projections $p'_j$ in $A_\omega$ by 
\[ p'_j=p_j+\sum_{l=0}^{L-1}\left(e_{lM+j}-\alpha_1^{lM+j}(ww^*)\right). \]
Then, we can check 
\[ \sum_{j=0}^{M-1}p'_j+\sum_{j=0}^Mq_j=1, \]
\[ \alpha_i(p'_j)=p_j \]
for all $i=2,3,\dots,N$ and $j=0,1,\dots,M-1$ and 
\[ \lVert\alpha_1(p'_j)-p'_{j+1}\rVert<3\ep+\frac{4}{\sqrt{K}} \]
for each $j=0,1,\dots,M-1$, where $p'_M$ is understood as $p'_0$. 
This completes the proof. 
\end{proof}

\section{Classification}

In this section, we would like to prove our main result Theorem \ref{main} 
by using the Evans-Kishimoto intertwining argument (\cite{EK}). 

\begin{lem}\label{1-cocycle}
Let $A$ be a unital $C^*$-algebra which is isomorphic to $\mathcal{O}_2$ 
and let $\alpha$ and $\beta$ be actions of $\Z^N$ on $\mathcal{O}_2$. 
When $\alpha$ has the Rohlin property, we have the following. 
\begin{enumerate}
\item There exists an $\alpha$-cocycle $\{u_n\}_n$ in $A^\omega$ 
such that $\beta_n(a)=\Ad u_n\circ\alpha_n(a)$ 
for all $n\in\Z^N$ and $a\in A$. 
\item For any $\ep>0$ and a finite subset $\mathcal{F}$ of $A$, 
there exists an $\alpha$-cocycle $\{u_n\}_n$ in $A$ 
such that 
\[ \lVert\beta_{\xi_i}(a)-\Ad u_{\xi_i}\circ\alpha_{\xi_i}(a)\rVert<\ep \]
for each $i=1,2,\dots,N$ and $a\in\mathcal{F}$. 
\end{enumerate}
\end{lem}
\begin{proof}
To simplify notation, 
we denote $\alpha_{\xi_i}$, $\beta_{\xi_i}$ by $\alpha_i$, $\beta_i$. 

(1). 
We prove this by induction on $N$. 
When $N=1$, the assertion is clearly true, 
because $\alpha_1,\beta_1\in\Aut(A)$ are approximately unitarily equivalent 
(\cite[Theorem 3.6]{R1}). 
Suppose that the assertion for $N-1$ has been shown. 
Let us consider the case of $\Z^N$-actions. 
Let $\alpha'$ be the action of $\Z^{N-1}$ on $A$ 
generated by $\alpha_2,\alpha_3,\dots,\alpha_N$. 
By using the induction hypothesis to $\alpha'$, 
we can find unitaries $u_2,u_3,\dots,u_N$ in $A^\omega$ such that 
\[ \beta_i(a)=\Ad u_i\circ\alpha_i(a) \]
and 
\[ u_i\alpha_i(u_j)=u_j\alpha_j(u_i) \]
for all $i,j=2,3,\dots,N$ and $a\in A$. 
Since two automorphisms $\alpha_1$ and $\beta_1$ on $A$ 
are approximately unitarily equivalent, 
there exists a unitary $u\in A^\omega$ such that 
\[ \beta_1(a)=\Ad u\circ\alpha_1(a) \]
for all $a\in A$. 
For $i=2,3,\dots,N$, we define $x_i\in A^\omega$ by 
\[ x_i=u\alpha_1(u_i)(u_i\alpha_i(u))^*. \]
It is easy to see that $x_i$ belongs to $\mathcal{U}(A_\omega)$. 
Let $\tilde\alpha$ be the perturbed action of $\alpha'$ 
by the $\alpha'$-cocycle $\{u_2,u_3,\dots,u_N\}$. 
Then, we can verify, for every $i,j=2,3,\dots,N$, 
\begin{align*}
x_i\tilde\alpha_i(x_j)
&=u\alpha_1(u_i)\alpha_i(u)^*u_i^*
\tilde\alpha_i(u\alpha_1(u_j)\alpha_j(u)^*u_j^*) \\
&=u\alpha_1(u_i)\alpha_i(\alpha_1(u_j)\alpha_j(u)^*u_j^*)u_i^* \\
&=u\alpha_1(u_i\alpha_i(u_j))\alpha_i(\alpha_j(u)^*)(u_i\alpha_i(u_j))^* \\
&=u\alpha_1(u_j\alpha_j(u_i))\alpha_j(\alpha_i(u)^*)(u_j\alpha_j(u_i))^* \\
&=x_j\tilde\alpha_j(x_i), 
\end{align*}
and so the family of unitaries $\{x_2,x_3,\dots,x_N\}$ is 
a $\tilde\alpha$-cocycle on $A_\omega$. 
It follows from Theorem \ref{CVanish} that 
there exists a unitary $v\in A_\omega$ such that 
\[ x_i=v\tilde\alpha_i(v^*) \]
for all $i=2,3,\dots,N$. 
Put $u_1=v^*u$. 
It can be easily checked that 
\[ \beta_1(a)=\Ad u_1\circ\alpha_1(a) \]
for all $a\in A$ and 
\[ u_1\alpha_1(u_i)=u_i\alpha_i(u_1) \]
for each $i=2,3,\dots,N$. 
Therefore, the family of unitaries $\{u_1,u_2,\dots,u_N\}$ 
induces the desired $\alpha$-cocycle in $A^\omega$. 

(2). 
From (1), 
there exists an $\alpha$-cocycle $\{u_n\}_n$ in $A^\omega$ 
such that 
\[ \beta_n(a)=\Ad u_n\circ\alpha_n(a) \]
for all $n\in\Z^N$ and $a\in A$. 
It follows from Theorem \ref{CVanish} that 
there exists a unitary $v\in A^\omega$ such that 
\[ u_n=v\alpha_n(v^*) \]
for every $n\in\Z^N$. 
Hence, for any $\ep>0$ and a finite subset $\mathcal{F}$ of $A$, 
we can find a unitary $w$ in $A$ such that 
\[ \lVert \beta_i(a)-\Ad(v\alpha_i(v^*))\circ\alpha_i(a)\rVert<\ep \]
for every $i=1,2,\dots,N$ and $a\in\mathcal{F}$. 
Therefore, $\{v\alpha_n(v^*)\}_n$ is the desired $\alpha$-cocycle. 
\end{proof}

Now we are ready to give a proof for our main theorem. 
We make use of the Evans-Kishimoto intertwining argument 
(\cite[Theorem 4.1]{EK}). 
See also \cite[Theorem 5]{N2} or \cite[Theorem 3.5]{I2}. 

\begin{thm}\label{main}
Let $\alpha$ and $\beta$ be two outer actions 
on the Cuntz algebra $\mathcal{O}_2$. 
Then, they are cocycle conjugate to each other. 
\end{thm}
\begin{proof}
We denote the Cuntz algebra $\mathcal{O}_2$ by $A$. 
Set $S=\{\xi_1,\xi_2,\dots,\xi_N\}\subset\Z^N$. 
Note that, by Theorem \ref{Rohlintype}, 
both $\alpha$ and $\beta$ have the Rohlin property. 
We choose an increasing family of finite subsets 
$\mathcal{F}_1,\mathcal{F}_2,\dots$ of $A$ 
whose union is dense in $A$. 
Put $\alpha^0=\alpha$ and $\beta^1=\beta$. 
We will construct 
$\Z^N$-actions $\alpha^{2k}$ and $\beta^{2k+1}$ on $A$, 
cocycles $\{u^k_n\}_{n\in\Z^N}$ in $A$ 
and unitaries $v_k$ in $A$ inductively. 

Applying Corollary \ref{appCVanish} 
to $\beta^1$, $2^{-1}>0$ and $\mathcal{F}_1$, 
we obtain $\delta_1>0$ and a finite subset $\mathcal{G}_1$ of $A$. 
We let 
\[ \mathcal{G}_1'=\bigcup_{g\in S}\beta^1_{-g}(\mathcal{G}_1)
\cup\mathcal{F}_1. \]
By Lemma \ref{1-cocycle} (2), 
there exists an $\alpha^0$-cocycle $\{u^0_n\}_n$ in $A$ 
such that 
\begin{equation}
\lVert\beta^1_g(a)-\Ad u^0_g\circ\alpha^0_g(a)\rVert
<\frac{\delta_1}{2}
\label{first}
\end{equation}
for every $g\in S$ and $a\in\mathcal{G}_1'$. 
Let $\alpha^2$ be the perturbed action 
of $\alpha^0$ by the $\alpha^0$-cocycle $\{u^0_n\}_n$. 
By using Corollary \ref{appCVanish} 
for $\alpha^0$ and $u^0$, 
we can find a unitary $v_0$ in $A$ such that 
\[ \lVert u^0_g-v_0\alpha^0_g(v_0)^*\rVert<1 \]
for each $g\in S$. 

Applying Corollary \ref{appCVanish} 
to $\alpha^2$, $2^{-2}$ and 
\[ \mathcal{F}_2'=\mathcal{F}_2\cup\Ad v_0(\mathcal{F}_2), \]
we obtain $\delta_2>0$ and a finite subset $\mathcal{G}_2$ of $A$. 
We may assume that $\delta_2$ is less than $\delta_1$ and $2^{-2}$. 
We let 
\[ \mathcal{G}_2'=\bigcup_{g\in S}\alpha^2_{-g}(\mathcal{G}_2)
\cup\bigcup_{g\in S}\beta^1_{-g}(\mathcal{G}_1)
\cup\mathcal{F}_2. \]
By Lemma \ref{1-cocycle} (2), 
there exists an $\beta^1$-cocycle $\{u^1_n\}_n$ in $A$ 
such that 
\begin{equation}
\lVert\Ad u^1_g\circ\beta^1_g(a)-\alpha^2_g(a)\rVert
<\frac{\delta_2}{2}
\label{second}
\end{equation}
for every $g\in S$ and $a\in\mathcal{G}_2'$. 
Let $\beta^3$ be the perturbed action 
of $\beta^1$ by the $\beta^1$-cocycle $\{u^1_n\}_n$. 
From \eqref{first} and \eqref{second}, one has 
\[ \lVert[u^1_g,a]\rVert<\delta_1 \]
for every $g\in S$ and $a\in\mathcal{G}_1$. 
By using Corollary \ref{appCVanish} 
for $\beta^1$ and $u^1$, 
we can find a unitary $v_1$ in $A$ such that 
\[ \lVert u^1_g-v_1\beta^1_g(v_1)^*\rVert<2^{-1} \]
and 
\[ \lVert[v_1,a]\rVert<2^{-1} \]
for each $g\in S$ and $a\in\mathcal{F}_1$. 

Applying Corollary \ref{appCVanish} 
to $\beta^3$, $2^{-3}>0$ and 
\[ \mathcal{F}_3'=\mathcal{F}_3\cup\Ad v_1(\mathcal{F}_3), \]
we obtain $\delta_3>0$ and a finite subset $\mathcal{G}_3$ of $A$. 
We may assume that $\delta_3$ is less than $\delta_2$ and $2^{-3}$. 
We let 
\[ \mathcal{G}_3'=\bigcup_{g\in S}\beta^3_{-g}(\mathcal{G}_3)
\cup\bigcup_{g\in S}\alpha^2_{-g}(\mathcal{G}_2)
\cup\mathcal{F}_3. \]
By Lemma \ref{1-cocycle} (2), 
there exists an $\alpha^2$-cocycle $\{u^2_n\}_n$ in $A$ 
such that 
\begin{equation}
\lVert\beta^3_g(a)-\Ad u^2_g\circ\alpha^2_g(a)\rVert
<\frac{\delta_3}{2}
\label{third}
\end{equation}
for every $g\in S$ and $a\in\mathcal{G}_3'$. 
Let $\alpha^4$ be the perturbed action 
of $\alpha^2$ by the $\alpha^2$-cocycle $\{u^2_n\}_n$. 
From \eqref{second} and \eqref{third}, one has 
\[ \lVert[u^2_g,a]\rVert<\delta_2 \]
for every $g\in S$ and $a\in\mathcal{G}_2$. 
By using Corollary \ref{appCVanish} 
for $\alpha^2$ and $u^2$, 
we can find a unitary $v_2$ in $A$ such that 
\[ \lVert u^2_g-v_2\alpha^2_g(v_2)^*\rVert<2^{-2} \]
and 
\[ \lVert[v_2,a]\rVert<2^{-2} \]
for each $g\in S$ and $a\in\mathcal{F}_2'$. 

Applying Corollary \ref{appCVanish} 
to $\alpha^4$, $2^{-4}$ and 
\[ \mathcal{F}_4'=\mathcal{F}_4
\cup\Ad(v_2v_0)(\mathcal{F}_4), \]
we obtain $\delta_4>0$ and a finite subset $\mathcal{G}_4$ of $A$. 
We may assume that $\delta_4$ is less than $\delta_3$ and $2^{-4}$. 
We let 
\[ \mathcal{G}_4'=\bigcup_{g\in S}\alpha^4_{-g}(\mathcal{G}_4)
\cup\bigcup_{g\in S}\beta^3_{-g}(\mathcal{G}_3)
\cup\mathcal{F}_4. \]
By Lemma \ref{1-cocycle} (2), 
there exists an $\beta^3$-cocycle $\{u^3_n\}_n$ in $A$ 
such that 
\begin{equation}
\lVert\Ad u^3_g\circ\beta^3_g(a)-\alpha^4_g(a)\rVert
<\frac{\delta_4}{2}
\label{fourth}
\end{equation}
for every $g\in S$ and $a\in\mathcal{G}_4'$. 
Let $\beta^5$ be the perturbed action 
of $\beta^3$ by the $\beta^3$-cocycle $\{u^3_n\}_n$. 
From \eqref{third} and \eqref{fourth}, one has 
\[ \lVert[u^3_g,a]\rVert<\delta_3 \]
for every $g\in S$ and $a\in\mathcal{G}_3$. 
By using Corollary \ref{appCVanish} 
for $\beta^3$ and $u^3$, 
we can find a unitary $v_3$ in $A$ such that 
\[ \lVert u^3_g-v_3\beta^3_g(v_3)^*\rVert<2^{-3} \]
and 
\[ \lVert[v_3,a]\rVert<2^{-3} \]
for each $g\in S$ and $a\in\mathcal{F}_3'$. 

Repeating this argument, 
we obtain a sequence of $\Z^N$-actions 
$\alpha^0,\beta^1,\alpha^2,\beta^3,\dots$, 
cocycles $\{u^0_n\}_n,\{u^1_n\}_n,\dots$ 
and unitaries $v_0,v_1,\dots$. 
Define $\sigma_{2k}$ and $\sigma_{2k+1}$ by 
\[ \sigma_{2k}=\Ad(v_{2k}v_{2k-2}\dots v_0) \]
and 
\[ \sigma_{2k+1}=\Ad(v_{2k+1}v_{2k-1}\dots v_1). \]
Since we have 
\[ \lVert[v_k,a]\rVert<2^{-k} \]
and 
\[ \lVert[v_k,\sigma_{k-2}(a)]\rVert<2^{-k} \]
for any $a\in\mathcal{F}_k$ , 
we can conclude that 
there exist automorphisms $\gamma_0$ and $\gamma_1$ such that 
\[ \gamma_0=\lim_{k\to\infty}\sigma_{2k} \]
and 
\[ \gamma_1=\lim_{k\to\infty}\sigma_{2k+1} \]
in the point-norm topology (see \cite[Lemma 3.4]{I2}). 

Define $w^{2k}_g,w^{2k+1}_g\in\mathcal{U}(A)$ by 
\[ w^{2k}_g=u^{2k}_g\alpha^{2k}_g(v_{2k})v_{2k}^* \]
and 
\[ w^{2k+1}_g=u^{2k+1}_g\beta^{2k+1}_g(v_{2k+1})v_{2k+1}^* \]
for every $k=0,1,2,\dots$ and $g\in S$. 
Then $\lVert w^k_g-1\rVert$ is less than $2^{-k}$. 
Furthermore, we define $\widetilde{w}^k_g\in\mathcal{U}(A)$ by 
\[ \widetilde{w}^0_g=w^0_g, \quad \widetilde{w}^1_g=w^1_g \]
and 
\[ \widetilde{w}^k_g=w^k_gv_k\widetilde{w}^{k-2}_gv_k^* \]
inductively. 
We would like to see that 
$\{\widetilde{w}^{2k}_g\}_k$ converges to a unitary for each $g\in S$. 
From 
\[ \widetilde{w}^{2k}_g=w^{2k}_g\cdot\Ad(v_{2k})(w^{2k-2}_g)
\cdot\Ad(v_{2k}v_{2k-2})(w^{2k-4}_g)\cdot\dots
\cdot\Ad(v_{2k}v_{2k-2}\dots v_2)(w^0_g), \]
we get 
\[ \sigma_{2k}^{-1}(\widetilde{w}^{2k}_g)
=\sigma_{2k-2}^{-1}(w^{2k-2}_g)
\cdot\sigma_{2k-4}^{-1}(w^{2k-4}_g)\cdot\dots
\cdot\sigma_0^{-1}(w^0_g). \]
It follows from $\lVert w^{2k}_g-1\rVert<4^{-k}$ that 
the right hand side converges. 
Hence $\widetilde{w}^{2k}_g$ converges to a unitary $W^0_g$ in $A$, 
because $\sigma_{2k}$ converges to $\gamma_0$. 
Likewise, 
$\widetilde{w}^{2k+1}_g$ also converges to a unitary $W^1_g$ in $A$. 
We can also check that 
the unitaries $\{\widetilde{w}^{2k}_g\}_g$ are a cocycle 
for the $\Z^N$-action $\sigma_{2k}\circ\alpha\circ\sigma_{2k}^{-1}$ 
and that 
the unitaries $\{\widetilde{w}^{2k+1}_g\}_g$ are a cocycle 
for the $\Z^N$-action $\sigma_{2k+1}\circ\beta\circ\sigma_{2k+1}^{-1}$. 
In addition, we can verify 
\[ \Ad(\widetilde{w}^{2k}_g)\circ
\sigma_{2k}\circ\alpha_g\circ\sigma_{2k}^{-1}=\alpha^{2k+2}_g \]
and 
\[ \Ad(\widetilde{w}^{2k+1}_g)\circ
\sigma_{2k+1}\circ\beta_g\circ\sigma_{2k+1}^{-1}=\beta^{2k+3}_g. \]
Since 
\[ \lVert\beta^{2k+3}_g(a)-\alpha^{2k+2}_g(a)\rVert<2^{-2k-3} \]
for all $a\in\mathcal{F}_{2k+2}$, 
we obtain 
\[ \Ad W^1_g\circ\gamma_1\circ\beta_g\circ\gamma_1^{-1}
=\Ad W^0_g\circ\gamma_0\circ\alpha_g\circ\gamma_0^{-1} \]
for every $g\in S$. 
Furthermore, $\{W^0_g\}_g$ is a cocycle 
for the $\Z^N$-action $\gamma_0\circ\alpha\circ\gamma_0^{-1}$ and 
$\{W^1_g\}_g$ is a cocycle 
for the $\Z^N$-action $\gamma_1\circ\beta\circ\gamma_1^{-1}$. 
Therefore, we can conclude that 
$\alpha$ and $\beta$ are cocycle conjugate to each other. 
\end{proof}

\begin{rem}
From Theorem \ref{main} and Theorem \ref{CVanish}, 
one can actually show the following: 
Let $\alpha$ and $\beta$ be outer actions of $\Z^N$ on $\mathcal{O}_2$. 
For any $\ep>0$, 
there exist $\gamma\in\Aut(\mathcal{O}_2)$ and 
an $\alpha$-cocycle $\{u_n\}_{n\in\Z^N}$ such that 
\[ \Ad u_n\circ\alpha_n(a)=\gamma\circ\beta_n\circ\gamma^{-1}(a), \]
\[ \lVert u_{\xi_i}-1\rVert<\ep \]
for any $n\in\Z^N$, $a\in A$ and $i=1,2,\dots,N$. 
\end{rem}

\end{document}